\documentclass[a4paper,english,cleveref]{lipics-v2021}
\hideLIPIcs
\nolinenumbers

\usepackage{amsthm}
\usepackage{amsmath}
\usepackage{amssymb}
\usepackage{amsfonts}
\usepackage{epsfig}
\usepackage{graphics}
\usepackage{graphicx}
\usepackage{booktabs}
\usepackage{pdfpages}
\usepackage{MnSymbol}

\usepackage{complexity}
\usepackage{url}

\usepackage{latexsym}
\usepackage{psfrag}
\usepackage{enumerate}
\usepackage{here}
\usepackage{esvect}
\usepackage[update,prepend]{epstopdf}
\usepackage{anyfontsize}
\usepackage{url}
\usepackage{hyperref}

\usepackage{microtype}%if unwanted, comment out or use option "draft"

\newcommand{\etal}{{et~al.}}
\newcommand{\ie}{{i.e.}}

\newcommand{\alg}{\textsf{ALG}}
\newcommand{\opt}{\textsf{OPT}}

\newcommand{\RR}{\mathbb{R}} %  set of real numbers
\newcommand{\eps}{\varepsilon}

\providecommand{\intd}[0]%
{\;\mbox{d}}

\usepackage[noline,linesnumbered]{algorithm2e}
\SetArgSty{textrm}
\let\oldnl\nl
\newcommand{\nonl}{\renewcommand{\nl}{\let\nl\oldnl}}
\makeatletter
\def\TitleOfAlgo{\@ifnextchar({\@TitleOfAlgoAndComment}{\@TitleOfAlgoNoComment}}
\def\@TitleOfAlgoAndComment(#1)#2{\nonl\hspace*{-1.5em}#2 #1\;}
\def\@TitleOfAlgoNoComment#1{\nonl\hspace*{-1.5em}#1\;}
\makeatother
\DontPrintSemicolon

\newcommand{\later}[1]{}
\newcommand{\old}[1]{}

\theoremstyle{} % theorems enumerated by letters
    \newtheorem{prevthm}{Theorem}

\title{Closed Curve Covering and \\ Multiagent TSP Ratios}

\titlerunning{Closed Curve Covering and Multiagent TSP Ratios}

\author{Travis Dillon}
{Massachusetts Institute of Technology}
{travis.dillon@mit.edu}
{0000-000?-4250-8836}
{}

\author{Adrian Dumitrescu}
{Algoresearch L.L.C., Milwaukee, WI, USA}
{ad.dumitrescu@algoresearch.org}
{0000-0002-1118-0321}
{}

\authorrunning{Travis Dillon and Adrian Dumitrescu}

\Copyright{Travis Dillon and Adrian Dumitrescu}

\funding{TD's work was partially supported by a National Science Foundation
Graduate Research Fellowship under Grant No. 2141064.}



\begin{document}

\maketitle

\begin{abstract}
  How efficiently can a closed curve of unit length in $\RR^d$ be covered by $k$ closed curves so as to minimize the
  maximum length of the $k$ curves? We show that the maximum length is at most $2k^{-1} -  \frac{1}{4} k^{-4}$
  for all $k\geq 2$ and $d \geq 2$.
  
  As a first byproduct, we show that $k$ agents can traverse a  Euclidean TSP instance significantly faster than a single agent.
  We thereby sharpen recent planar results by Berendsohn, Kim, and Kozma (2025) and extend these improvements to all dimensions.

  As a second byproduct, we obtain a linear time approximation algorithm with ratio $2 - \frac{1}{4} k^{-3}$ for covering
  any closed polygonal curve in $\RR^d$ by $k$ closed curves so that the maximum length of an individual curve is minimized.  

\end{abstract}

\section{Introduction}\label{sec:intro}

Many problems from the manufacturing industry involve dividing a task into a number of smaller tasks
while optimizing a specific objective function.
As one example, consider the problem of partitioning a given set $P$
of points in the plane into $k$  subsets, $P_1,\ldots,P_k$, such that
$\sum_{i=1}^k  t(P_i)$ is minimized, where $t(\cdot)$ is the objective function. For example,
$t(\cdot)$ may be the length of a \emph{minimum spanning tree} (MST) 
or the length of  a  shortest traveling salesman tour (TSP tour) of the point set.
Alternatively, one may want to minimize $\max_{i=1}^k  t(P_i)$.
Andersson, Gudmundsson, Levcopoulos, and Narasimhan~\cite{AGLN03} studied the sum variant for the MST,
and mentioned that variants of this problem arise in applications from the shipbuilding industry.

Here we investigate a related problem: Given a closed curve in $\RR^d$, the goal is
to cover it by $k$ closed curves so as to minimize the sum, or the maximum, of the lengths of the $k$ curves.
For any closed curve $C$, let $P_k(C)$ denote the set of partitions of $C$ into $k$ (not necessarily connected) non-empty sets,
while $\hat P_k(C)$ denotes the partitions into equal-length sets.\footnote{The reason we require a partition into equal-length
  sets is that without this restriction the minimum value of $\beta(\mathcal S)$ is obtained when one part of the partition is $P$
  and the other parts have length 0.}
For any $\mathcal S \in P_k(C)$,  we define 
\[
    \beta(\mathcal S) = \frac{\textstyle \sum_{S \in \mathcal S} w(S)}{w(C)} \text{~~~~and~~~~}
    \beta(C,k) = \inf_{\mathcal S \in \hat P_k(C)} \beta(\mathcal S).
\]

If $C$ is a long, thin rectangle, then $\beta(C,k)$ can be arbitrarily close to 1 (and $\beta(\mathcal S) \geq 1$ always).
Our first result is a precise determination of the upper bound for $\beta(C,k)$ among all possible curves $C$.
Let 
\[
    \beta_d(k) := \sup_{C \text{ in } \RR^d} \beta(C,k).
\]

\begin{theorem}\label{thm:max-avg-tsp}
    $\beta_d(k) = \frac{1}{k} + \frac{1}{\pi} \sin \frac{\pi}{k}$ for every $k \geq 1$ and $d \geq 2$.
\end{theorem}

Recently, Berendsohn, Kim, and Kozma \cite{BKK25} proved bounds on the
ratio between the \emph{longest} of the $k$ smaller cycles and $C$. Their results are phrased in terms of point sets
and minimal closed curves containing them. (In terms of the first paragraph of this section, they take $t(P)$ to be
the length of the smallest closed curve containing $P$.) This perspective makes a close
connection with the Euclidean Traveling Salesman Problem, which has been the
subject of intense study for many years~\cite{Ar98,BHH59,GGJ76,Mi99,Pa77}.
Their analysis focuses on partitioning the minimal closed curve itself, rather than the point set;
since any closed curve is arbitrarily close to the minimal cycle through some finite point set,
we lose no generality by focusing exclusively on curves.

In essence, they defined 
\[
    \gamma(C,k)
    := \min_{\mathcal S \in P_k(C)}  \frac{\max_{S \in \mathcal S} w(S)}{w(C)},
\]
and obtained several results on the extremal ratio
\[
    \gamma_d(k) := \sup_{C \text{ in } \RR^d} \gamma(C,k)
\]
when $d=2$, including the exact value for $\gamma_2(2)$.

\begin{figure}[htbp]
 \centering
 \includegraphics[scale=0.83]{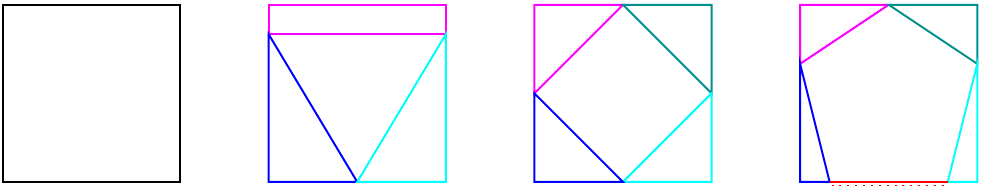}
 \caption{When $C$ is the unit square, we have: $\gamma(C,3) \leq 0.580$, with the splitting points $(0.5,0)$, $(0,y_0)$, and $(1,y_0)$,
   where $y_0>1/2$ is a solution to $8y^2 -21 y +12=0$; 
   $\gamma(C,4) \leq (2+\sqrt{2})/8 \leq 0.427$, with the splitting points as the midpoints of the four sides;
   and $\gamma(C,5) \leq 0.366$ (details omitted).}
 % xxx: need to check these values /upper bounds 
 \label{fig:square}
\end{figure}

A few upper bounds on $\gamma(C,k)$, \ie, covering $C$ by closed curves when $C$ is a square and $k=3,4,5$,
are shown in Fig.~\ref{fig:square}.

\begin{prevthm}[Berendsohn, Kim, Kozma \cite{BKK25}]\label{thm:BKK-bound}
    $\gamma_2(k) \geq \frac{1}{k} + \frac{1}{\pi} \sin \frac{\pi}{k}$.
  When $k=2$, this is tight: $\gamma_2(2) = \frac{1}{2} + \frac{1}{\pi}$.
\end{prevthm}

The lower bound comes from a circle of unit length, in which case the optimal partition consists
of equal arcs of length $1/k$. 

They also obtained upper bounds for $k \geq 2$, by iteratively cutting the closed curve by chords:
\begin{itemize}
    \item $\gamma_2(a\cdot b) \leq \gamma_2(a) \gamma_2(b)$
    \item $\gamma_2(a+b) \leq \big(1 + \frac{2}{\pi}\big)
      \displaystyle\frac{\gamma_2(a) \gamma_2(b)}{\gamma_2(a) + \gamma_2(b)}.$
\end{itemize}

Using these formulas, they obtain the bounds in \cref{table:bounds}.
The upper and lower bounds appear to grow increasingly farther apart:
by $k=10$, the bounds differ by a factor of more than 2.5.
Indeed, iterated application of either upper bound yields
$\gamma_2(2^n) \leq \gamma_2(2)^n \approx (0.818)^n$, while the lower bound
is $\gamma_2(2^n) \geq \frac{1}{2^n} + \frac{1}{\pi} \sin \frac{\pi}{2^n} \approx 2^{-n+1}$.
Under the plausible assumption that the second bound is optimized when $a=b$
(which is how the upper bounds for all even $k \leq 10$ are obtained),
the ratio between the lower and upper bounds for $\gamma_2(k)$ grows polynomially with $k$.  

\begin{table}[htbp]
\centering
\scalebox{0.93}{
    \begin{tabular}{|c|c c c c c c c c c c|}\hline
    $k$ & 1 & 2 & 3 & 4 & 5 & 6 & 7 & 8 & 9 & 10 \\[0.5ex] \hline
    Lower bound \cite{BKK25} & 1 & $0.818$ & $0.609$ & $0.475$ & $0.387$ & $0.325$ & $0.281$ & $0.246$ & $0.220$ & $0.198$ \\
    Upper bound \cite{BKK25} ($d=2$) & 1 & $0.818$ & $0.737$ & $0.670$ & $0.634$ & $0.603$ & $0.574$ & $0.548$ & $0.533$ & $0.519$ \\
    New upper bound  & -- & -- & $0.644$ & $0.494$ & $0.398$ & $0.333$ & $0.285$ & $0.250$ & $0.222$ & $0.200$\\\hline
      \end{tabular}}
\vspace{0.1in}
\caption{Lower and upper bounds for $\gamma_d(k)$ from this paper and~\cite{BKK25}.
  Upper bounds are rounded up; lower bounds are rounded down.}
\label{table:bounds}
\end{table}

In \cite{BKK25}, the authors propose further investigation of $\gamma_2(k)$ for $k\geq 3$ as $\gamma_d(k)$ for $d \geq 3$.
The starting point for our results is a simple, asymptotically optimal bound for $\gamma_d(k)$ in all dimensions,
which moreover substantially improves the upper bounds for $\gamma_2(k)$ for all $k\geq 3$.

\begin{proposition}\label{prop:2/k}
    $\gamma_d(k) \leq 2/k$ for all $k\geq 2$ and $d \geq 2$.
\end{proposition}

Because $\RR^2$ embeds in $\RR^d$, \cref{thm:BKK-bound} implies that
$\gamma_d(k) \geq \gamma_2(k) \geq \frac{1}{k} + \frac{1}{\pi} \sin \frac{\pi}{k}$,
which tends to $2/k$ as $k \to \infty$. Therefore $\gamma_d(k) \sim 2/k$ as
$k\to\infty$ for every $d$. The remaining additive gap between the lower and
upper bounds for $\gamma_d(k)$ is quite small, namely $\Theta(k^{-3})$. Still,
it would be interesting to know the exact value: is it
$\frac{1}{k} + \frac{1}{\pi}\sin \frac{\pi}{k}$? 

While we do not resolve this question, by using a slightly non-uniform partition of the curve,
we do show that the upper bound in Proposition~\ref{prop:2/k} can be improved.
Our second result is the following. 

\begin{theorem}\label{thm:<2/k}
  $\gamma_d(k) \leq \frac{2}{k}-  \frac{1}{4k^4}$ for all $k\geq 3$ and $d \geq 2$.
  Furthermore, for $k=3,\ldots,10$, $\gamma_d(k)$ is bounded from above as in row $3$ of Table~\ref{table:bounds}.
\end{theorem}

The new upper bounds $3 \leq k \leq 10$ are obtained by numerically optimizing our argument
for those specific values of $k$.

We then turn to the problem of computing a $k$-covering of a closed polygonal curve:
Given a closed polygonal curve with $n$ edges, compute $k$ closed polygonal curves that cover it,
and such that the maximum individual length is minimized.

The upper bounds from Berendsohn, Kim, and Kozma, since they're recursive, could be used to implement
a recursive algorithm for this problem. This algorithm, at each step, would need to find a direction
in which the curve has small width. Our $2/k$ result gives an almost trivial algorithm: for $i=0,\ldots,k-1$,
take the polygonal arcs between $x = i/k$ and $x=(i+1)/k$, and make them into $k$ closed polygonal curves.
This is faster, of course, and with better approximation ratio: indeed, we have $\opt \geq 1/k$, and $\alg \leq 2/k$
by the triangle inequality, thus the approximation ratio $2$ is immediate.

We further improve the ratio $2$ to $2 - \frac{1}{4} k^{-3}$, based on Theorem~\ref{thm:<2/k}.
For $k \leq 10$, the approximation ratios are listed below.

\begin{table}[htbp]
\centering
\scalebox{0.95}{
    \begin{tabular}{|c|c c c c c c c c c c|}\hline
    $k$ & 1 & 2 & 3 & 4 & 5 & 6 & 7 & 8 & 9 & 10 \\[0.5ex] \hline
    Approximation ratio & $1$ & $1.637$ & $1.932$ & $1.974$ & $1.988$ & $1.993$ & $1.996$ & $1.998$ & $1.998$ & $2$ \\\hline
\end{tabular}}
% xxx: need to check these ratios
\vspace{0.1in}
\caption{Approximation ratios for $k=1,\ldots,10$; third digit rounded up.}
\label{table:ratios}
\end{table}

\begin{theorem}\label{thm:algo}
Given a closed polygonal chain $C=(p_1,p_2,\ldots,p_n)$ in $\RR^d$ and a positive integer $k \geq 3$,
a $2- \frac{1}{4k^3}$ approximation for covering $C$ by $k$ closed curves can be computed in $O(n+k)$ time.
Furthermore, for $k=3,\ldots,10$, the approximation ratios listed in Table~\ref{table:ratios} can be obtained.
\end{theorem}

\subparagraph{Related work.}
Estimating the length of a shortest tour of $n$ points in the unit square with respect to  Euclidean distances
has been studied as early as the 1940s by Fejes T\'oth~\cite{Fej40}, Few~\cite{Few55},
and Verblunsky~\cite{Ve51}, respectively.
Few~\cite{Few55} proved that the (Euclidean) length of a shortest cycle (tour)
through $n$ points in the unit square is at most $\sqrt{2n}+7/4$.
The same upper bound holds for the minimum spanning tree~\cite{Few55}.
Karloff~\cite{Ka89} derived a slightly better upper bound for the shortest cycle, $1.392 \sqrt{n} + 7/4$;
he also emphasized the difficulty of obtaining further improvements.

In the \emph{fence patrolling} problem, 
a set of $k$ mobile agents with maximum speeds
$v_i$ ($i=1,\ldots,k$) are in charge of patrolling a region of interest.
Czyzowicz~\etal~\cite{CGKK11} introduced a one-dimensional variant, where the agents move along a
rectifiable Jordan curve representing a \emph{fence}.
The movement of the agents is described by a \emph{patrolling schedule},
where the speed of the $i$th agent may vary between $0$ and its maximum value $v_i$ in any of the two directions along the fence.

Given a fence of length $1$ and maximum speeds $v_1,\ldots,v_k>0$, the problem is to find a
patrolling schedule that minimizes the \emph{idle time} $I$,
defined as the longest time interval in $[0,\infty)$ during which a
  point on the fence remains unvisited.
  A  volume argument~\cite{CGKK11} yields the lower bound
  $I \geq 1/\sum_{i=1}^{k} v_i$, which is $1/k$ for unit speeds.
An idle time of $2/k$ is always achievable with unit speeds.
This value coincides with the upper bound in Proposition~\ref{prop:2/k}. 
Other results on fence patrolling by multiple agents can be found in~\cite{Ch04,DT14,HKM+19,KK15}.

\section{Proofs of Theorems~\ref{thm:max-avg-tsp} and Theorems~\ref{thm:<2/k}} \label{sec:proofs}

\subparagraph{Proof of \cref{thm:max-avg-tsp}.}
Let $C$ be the unit circle. Berendsohn, Kim, and Kozma prove that $\gamma(C,k) \geq \frac{1}{k} + \frac{1}{\pi} \sin \frac{\pi}{k}$, and their proof also shows the same bound on $\beta(C,k)$. We produce here a different view of that proof, in terms of curves rather than points.

Let $\mathcal Q\in \hat P_k(C)$.
The minimum-length closed curve containing $Q_i \in \mathcal Q$
consists of a subset of the unit circle of length $2\pi/k$ together with several
chords of the circle. The chords span several arcs whose lengths $\alpha_1, \dots, \alpha_m$
sum to at least $2\pi/k$ (otherwise the length of $Q$ would be
less than $2\pi/k$). The length of a chord spanning an angle of
$\alpha \leq \pi$ is $2\sin(\alpha/2)$; because $\sin(\alpha/2) \geq \alpha/\pi \sin(\pi/k)$ for $0 \leq \alpha \leq 2\pi/k$,
the total length of the chords in $Q$ is at least $\sum_{i=1}^m 2\sin(\alpha_i/2) \geq 2\sin(\pi/k)$. 
 
In other words, $\sum_{Q \in \mathcal Q} w(Q) \geq 2\pi + 2k\sin(\pi/k)$ for any $\mathcal Q \in \hat P_k(C)$. Therefore
\[
    \beta(C,k)
    \geq \frac{2\pi + 2k\sin(\pi/k)}{2\pi k}
    = \frac{1}{k} + \frac{1}{\pi} \sin\frac{\pi}{k},
\]
which shows that $\beta_d(k) \geq \frac{1}{k} + \frac{1}{\pi} \sin \frac{\pi}{k}$.

For the upper bound, we partition $C$ into $k$ arcs of length $1/k$ with short chords.
Let $C$ be a closed curve of unit length in $\RR^d$, and let $r\colon [0,1] \to \RR^d$
be an arc-length parameterization of $C$. For each
$t \in [0,1]$ and $i \in \{0,1,\dots,k-1\}$, let $Q_{i,t}$ denote the arc of $r$ between
$x_i(t) := r(t+i/k)$ and $x_{i+1}(t) := r\big(t+(i+1)/k\big)$, and let $\mathcal Q_t = \{Q_{i,t}\}_{i=1}^k$. (The sums are
taken modulo 1.) Joining $x_i(t)$ and $x_{i+1}(t)$ by a straight segment makes a
closed curve containing $Q_{i,t}$, so $w(Q_{i,t}) \leq 1/k + \|x_i(t) - x_{i+1}(t)\|$.

We now pick $t$ uniformly at random in $[0,1]$. By linearity of expectation, we have
\begin{align*}
    \mathbb{E}\bigg[\frac{1}{k} \sum_{Q \in \mathcal Q_s} w(Q) \bigg]
    &\leq \frac{1}{k}\, \mathbb{E}\bigg[\sum_{Q \in \mathcal Q_s} \frac{1}{k} + \|x_i(s) - x_{i+1}(s)\|\bigg]\\
    &= \frac{1}{k} \bigg(1 + k\,\mathbb{E}\big[ \|x_0(s) - x_1(s)\|\big]\bigg).
\end{align*}
We now use the fact that a circle maximizes the average length of a chord spanning an arc of length $s$.

\begin{prevthm}[Abrams, Cantarella, Fu, Ghomi, and Howard \cite{ACFG03}]   \label{thm:avg}
For any $s \in [0,1/2]$ and any arc-length parameterization $r\colon [0,1] \to \RR^d$ of a length-1 closed curve, 
    \[
        \int_0^1 \|r(t+s) - r(t)\|\, dt \leq \frac{\sin{\pi s}}{\pi}.
    \]
    Equality holds if and only if $r$ parameterizes a circle.
\end{prevthm}

Exner, Harell, and Loss \cite{EHL06} independently obtained this result
a few years later. Both papers employ an argument that parallels the
Fourier analytic proof of the planar isoperimetric inequality. In any
case, the above result implies that 
\[
    \mathbb{E}\bigg[\frac{1}{k} \sum_{Q \in \mathcal Q_s} w(Q) \bigg]
    \leq \frac{1}{k} + \frac{1}{\pi} \sin\frac{\pi}{k}.
\]
In particular, there is a specific choice of $t \in [0,1]$ for which
$\frac{1}{k} \sum_{Q \in \mathcal Q_t} w(Q) \leq \frac{1}{k} + \frac{1}{\pi} \sin \frac{\pi}{k}$.
Thus, $\beta(C,k) \leq \frac{1}{k} + \frac{1}{\pi} \sin \frac{\pi}{k}$;
since $C$ was arbitrary, this proves \cref{thm:max-avg-tsp}.  
\qed

\subparagraph{Proof of Proposition~\ref{prop:2/k}.}
    Let $C$ be any closed curve of length 1, and let $\mathcal Q$ be a partition
    of $C$ into $k$ connected arcs of $C$, each of length $1/k$. Adding the
    segment joining the endpoints of $Q \in \mathcal Q$ forms each arc into a
    closed curve of length at most $2/k$ (by the triangle inequality), so $\gamma(C,k) \leq 2/k$. 
\qed

\subparagraph{Proof of Theorem~\ref{thm:<2/k}.}
Let $C$ be any closed curve of length 1.
Let
\[ s= \frac{1}{k} + (k-1) \eps, \text{ where } \eps= \frac{1}{8k^4}. \]
Our strategy is to take one arc of length $s$ and split the remaining portion of the curve into $k-1$ equal-length arcs.

By the series expansion of the $\sin x$ function around $x=0$, we have
\[ \sin x = x - \frac{x^3}{6} + \frac{x^5}{120} \cdots \leq   x - \frac{x^3}{12}, \text{ for } x \leq \pi. \]
By Theorem~\ref{thm:avg},
there exists an arc $a_1$ of $C$ of length $s$ which spans a chord whose length is at most
\begin{equation} \label{eq:a_1}
\frac{\sin{\pi s}}{\pi} \leq  \frac {\pi s}{\pi} - \frac {\pi^3 s^3}{ 12\pi}
= \left( 1 - \frac {\pi^2 s^2}{12} \right) s \leq  \left( 1 - \frac{1}{2k^2} \right) s.
\end{equation}
(We need $k \geq 3$ so that $s \leq 1/2$.) Consider the $k$-partition of $C$ into the arc $a_1$
and $k-1$ arcs $a_2,\ldots,a_k$ of the same length
    \[ \frac{1 - s}{k-1}= \frac{1}{k} -\eps. \]
We make each arc $a_i$ ($i = 1,\dots,k$) into a closed curve $C_i$ by adding the segment joining its endpoints.

The length of $C_1$ is at most
    \[ s + \left( 1 - \frac{1}{2k^2} \right) s
    = \frac{2}{k} + 2(k-1) \eps - \frac{1}{2k^3} - \frac{(k-1)\epsilon}{2k^2}
    \leq \frac{2}{k} - \frac{1}{4k^3}.   \]
    For $i \geq 2$, the length of $C_i$ is at most
   \[ \frac{2}{k} - 2\eps = \frac{2}{k} - \frac{1} {4k^4} \]
   by the triangle inequality. The length of each of the $k$ curves is therefore bounded above by
    \[ \max \left \{   \frac{2}{k} - \frac{1}{4k^3},\ \frac{2}{k} - \frac{1} {4k^4} \right \}
    = \frac{2}{k} - \frac{1} {4k^4}, \]
    as claimed.

To obtain specific upper bounds for $k=3,\ldots,10$, let $s_k$ be the solution of the equation
\begin{equation} \label{eq:sin}
  s + \frac{\sin{\pi s}}{\pi} = \frac{2(1-s)}{k-1}.
\end{equation}
    Following the previous argument with $s_k$ in place of $s$, we find a covering of $C$ with $k$ closed curves,
    each of length $2(1-s_k)/(k-1)$. These are the values that appear in \cref{table:bounds}.
\qed

\section{Approximation Algorithm: Proof of Theorem~\ref{thm:algo}}  \label{sec:algo}

In this section we consider closed polygonal curves with regard to the $k$-covering problem
and prove Theorem~\ref{thm:algo}.

\begin{proof}[Proof of Theorem \ref{thm:algo}]
We may assume that $C$ has unit length. Set $s = k^{-1} + (k-1) (8k^4)^{-1}$, 
%\[ a= \frac{1}{k} + \frac{k-1}{8k^4}, \]
as in the proof of Theorem~\ref{thm:<2/k}, and consider a moving interval of length $s$ on $C$,
initially with the left endpoint at a vertex of $C$. 

In the first phase, the algorithm computes an arc $a_1$ satisfying the inequality~\eqref{eq:a_1} by 
repeatedly advancing the moving interval counterclockwise on $C$ stopping at each of the $2n$ events when passing a
segment endpoint (until $C$ is completely traversed). We therefore obtain a set of at most $2n$ pairs of equal-length
intervals on $C$, where for each pair, the endpoints of the spanned chord belong to the same pair of edges,
say $e$ and $e'$, of the chain. The set of interval pairs can be computed in $O(n)$ time by a linear scan of $C$.
For each such pair $I=[a,b], I'=[a',b']$ of intervals on $C$, compute the corresponding shortest chord
and retain the overall shortest one in the end.

\begin{figure}[htbp]
 \centering
 \includegraphics[scale=0.8]{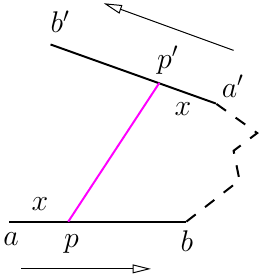}
 \caption{Minimizing the length of a moving chord spanning an arc of a fixed length on $C$.}
 \label{fig:chord}
\end{figure}

Computing the shortest chord for a given interval pair $I,I'$, amounts to minimizing a trigonometric function
of one variable (the position $x$ in $I$) and that depends solely on the pair of supporting edges $e,e'$.
Specifically, the algorithm computes the solution to the following optimization problem:
Given a sequence of four points $a,b,a',b'$ on $C$, where $ab \subset e$, $a'b' \subset e'$,
$|ab|=|a'b'|$, find $p \in ab$ and $p' \in a'b'$, such that $|ap|=|a'p'|=x$ (for some unknown $x$),
and $|p p'|$ is minimized. See Fig.~\ref{fig:chord}. Since $p$ and $p'$ move on two lines, the chord length
$|pp'|=|pp'(x)|$ is a \emph{unimodal} function with a unique minimum for $x \in [a,b]$.

A solution involves the following steps:
\begin{enumerate} \itemsep 1pt
\item Compute $|pa'|$ via the Law of Cosines in $\Delta{pba'}$.
\item Compute $\angle{pa'p'}$ via the  Law of Sines in $\Delta{pba'}$.
\item  Compute $|pp'|$ via the Law of Cosines in $\Delta{pa'p'}$.
\item Minimize $|pp'|$ as a function of $x$ (check the chords $aa'$ and $bb'$ against
  the one  given by the zero of the derivative). 
\end{enumerate}
Therefore, computing the solution for a given pair $I,I'$ of intervals takes $O(1)$ time.
 
Since there are at most $2n$ intervals, the first phase takes $O(n)$ time.
In the second phase, the algorithm subdivides the remaining part of the curve into $k-1$ chains of equal length;
this takes $O(n+k)$ time. Overall, the total time complexity is $O(n+k)$, as claimed.
(For $k=3,\ldots,10$, $s$ may be set according to~\eqref{eq:sin} if it yields an improvement.) 
\end{proof}

\begin{figure}[htbp]
 \centering
 \includegraphics[scale=0.95]{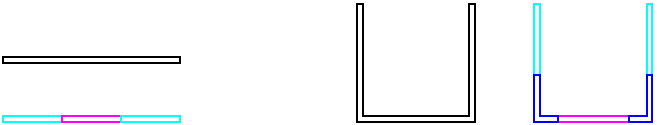}
 \caption{The hardness of computing an efficient curve covering.
   Left: a thin axis-parallel rectangle; here $k=3$.
   Right: a thin $U$-shaped axis-parallel polygonal curve; here $k=5$.}
 \label{fig:curves}
\end{figure}

\smallskip
Whereas we suspect that the $k$-covering problem for closed polygonal curves is $\NP$-hard, even for curves that
are  axis-parallel, it is possible that the problem admits better approximations for polygonal curves.
However, any algorithm that computes a set of $k$ covering curves by partitioning $C$ into
contiguous (or even connected) parts is doomed to avoid the vicinity of an optimal solution in some instances.
Let $\eps>0$ be small. Consider two axis-parallel curves $C_i$, $i=1,2$, shown in Fig.~\ref{fig:curves},
where the short sides have length $\eps$. Let $L_i$, $i=1,2$, denote the length of each curve.

(i) $C_1$ is a thin axis-parallel rectangle with two sides of unit length; so $L_1=2+2\eps$.
Let $k=3$. Whereas the covering sketched in  Fig.~\ref{fig:curves}\,(left) achieves
$\gamma(C_1,3) \leq 1/3 + O(\eps)$,
any covering obtained by partitioning $C_1$ into three contiguous parts only achieves
$\gamma'(C_1,3) \geq 1/2- O(\eps)$, where $\gamma'(C,k)$ stands for $\gamma(C,k)$ under this restriction.
(Note, this inequality is not the best possible, it just involves a simpler argument.)

To see this, one may assume that there is at most one splitting point on the upper side, 
thus one contiguous part is at least one half of one side, namely $1/2$, and so
the closed curve containing it has length at least $1$, and the claim follows.   

\smallskip
(ii) $C_2$ is a thin $U$-shape with three sides of unit length, two sides of length $1-\eps$
and one side of length $1-2\eps$; so $L_2=6 - 4\eps$.
Let $k=5$. Whereas the covering sketched in  Fig.~\ref{fig:curves}\,(right) achieves
$\gamma(C_2,5) \leq 1/5 + O(\eps)$,
any covering obtained by partitioning $C_2$ into five contiguous parts by five splitting points,
only achieves $\gamma'(C_2,5) \geq 1/3- O(\eps)$. To see this, partition $C_2$ into eight half-closed intervals,
by the vertices, going counterclockwise around $C_2$ (six long intervals, and two short intervals).
Since the number of splitting points is less than the number of long intervals, there exists one
empty long interval, and so the closed curve containing it has length at least $2-4\eps$,
and the claim follows.

\end{document}